\documentclass{ifacconf}

\usepackage{graphicx}
\usepackage{amsmath}
\usepackage{amssymb}
\usepackage{mathrsfs}
\usepackage{epstopdf}
\usepackage{caption}
\usepackage{subcaption}
\usepackage{balance}
\usepackage{natbib}

\newtheorem{theorem}{Theorem}

\newtheorem{remark}{Remark}
\newtheorem{proof}{Proof}

\begin{document}
\begin{frontmatter}

\title{Output-Feedback Stabilization of the Korteweg-de Vries Equation\thanksref{footnoteinfo}}

\thanks[footnoteinfo]{Financial support from Statoil ASA and the Norwegian Research Council (NFR project 210432/E30 Intelligent Drilling) is gratefully acknowledged. Corresponding email: \tt\small agusisma@itk.ntnu.no}

\author[First]{Agus Hasan}

\address[First]{Department of Engineering Cybernetics\\Norwegian University of Science and Technology\\Trondheim, Norway}

\begin{abstract}
The present paper develops boundary output-feedback stabilization of the Korteweg-de Vries (KdV) equation with sensors and an actuator located at different boundaries (anti collocated set-up) using backstepping method. The feedback control law and output injection gains are found using the backstepping method for linear KdV equation. The proof of stability is based on construction of a strict Lyapunov functional which includes the observer states. A numerical simulation is presented to validate the result.
\end{abstract}


\end{frontmatter}

\section{Introduction}

The Kortweg-de Vries (KdV) equation is a third-order partial differential equation (PDEs), which can be used to described weakly nonlinear shallow water surface \citep{KdV}. The KdV equation can be classified as hyperbolic-type PDEs, which describes a reversible dynamical process. Furthermore, it was found to have solitary wave solutions. The KdV equation is completely integrable and has infinitely many conserved quantities.

Boundary control of the KdV equation can be found in many literatures, e.g, \cite{Lionel2,Lionel1,Emman,KdV1,KdV2,KdV3}. In these literatures, the boundary control laws were found using the Lyapunov methods. Furthermore, only state-feedback was considered. Recent control design for the KdV equation includes the backstepping method \citep{KdV4,KdV5}, where the analysis mostly done for the linear KdV equation with state-feedback. In the infinite-dimensional backstepping, a Volterra integral transformation is used to transform the original system into a stable target system. Different with other approaches that require the solution of operator Riccati equations e.g., optimal control method \citep{Dic,AHasn}, backstepping yields control gain formulas which can be computed using symbolic computation and, in some cases, can even be given explicitly in terms of Bessel function \citep{kris} and Marcum Q-function \citep{Marcum}.

The backstepping method has been successfully used for control design of many PDEs such as the Schrodinger equation \citep{Guo}, the Ginzburg-Landau equation \citep{Aa}, the Navier-Stokes equation \citep{Vazque}, the surface wave equation \citep{AgusWave}, and the hyperbolic equation \citep{AgusPDE}. Furthermore, the backstepping method has found several applications in oil well drilling problems, including slugging control \citep{Flo1}, the lost circulation and kick problem \citep{Agus1,AHasan}, and the heave problem \citep{aamo,AgusDis}.

In this paper, we concern with the problem of output-feedback stabilization of the KdV equation using the infinite-dimensional backstepping method with sensors and an actuator located at different boundaries (anti collocated set-up). The design utilized the result for the linear KdV equation by \cite{Cerpa1}. The control law from the linear KdV is used to stabilize the systems, where the state is generated from a nonlinear observer of the KdV equation. The novelty of of this paper lies on the introduction of a strict Lyapunov functional for the output-feedback control problem. To prove the stability, we introduce a strict Lyapunov functional which equivalent to the $\mathbb{H}^3$ norm.

This paper is organized as follow. In section 2, we state the problem. Notations and definitions are presented in section 3. Output-feedback control for the linear KdV equation, which was solved by \cite{Cerpa1}, is presented in section 4. The main result for output-feedback stabilization of the KdV equation is presented in section 5. In section 6, we present a numerical example and the last section contains conclusions.

\section{Problem Statements}

We consider output-feedback stabilization of the Kortweg-de Vries (KdV) equation
\begin{eqnarray}
u_t(x,t)+u_x(x,t)+u_{xxx}(x,t)+u(x,t)u_x(x,t)&=&0\label{main}
\end{eqnarray}
with boundary conditions
\begin{eqnarray}
u(0,t)&=&U(t),\;\;u(1,t)=0,\;\;u_x(1,t)=0\label{bcmain3}
\end{eqnarray}
where $u:[0,1]\times[0,\infty)\rightarrow\mathbb{R}$. The subscripts $x$ and $t$ denote partial derivatives with respect to $x$ and $t$, respectively. The objective is to find a feedback law $U(t)$ to make the origin of \eqref{main}-\eqref{bcmain3} locally exponentially stable, using only measurements of $u_{xx}(1,t)$.

\section{Notations and Definitions}

For $u(x,t)\in\mathbb{R}$, we define
\begin{eqnarray}
\|u\|_{\infty} &=& \sup_{x\in[0,1]}|\gamma(x,t)|\\
\|u\|_{\mathbb{L}^1} &=& \int_0^1\; |\gamma(x,t)| \,\mathrm{d}x\\
\|u\|_{\mathbb{H}^i} &=& \int_0^1\; \sum_{j=0}^i\frac{\partial^j}{\partial x^j}\gamma(x,t) \,\mathrm{d}x
\end{eqnarray}
Furthermore, to simplify our notation, we denote $|u|=|u(x,t)|$ and $\|u\|=\|u(\cdot,t)\|$. For $u\in\mathbb{H}^3([0,1])$, recall the following well-known inequalities
\begin{eqnarray}
\|u\|_{\mathbb{L}^1} &\leq& c_1\|u\|_{\mathbb{L}^2} \leq c_2\|u\|_{\infty},\label{ineq1}\\
\|u\|_{\infty} &\leq& c_3\left(\|u\|_{\mathbb{L}^2}+\|u_x\|_{\mathbb{L}^2}\right) \leq c_4\|u\|_{\mathbb{H}^1},\label{ineq2}\\
\|u_x\|_{\infty} &\leq& c_5\left(\|u_x\|_{\mathbb{L}^2}+\|u_{xx}\|_{\mathbb{L}^2}\right) \leq c_6\|u\|_{\mathbb{H}^2},\label{ineq3}\\
\|u_{xx}\|_{\infty} &\leq& c_7\left(\|u_{xx}\|_{\mathbb{L}^2}+\|u_{xxx}\|_{\mathbb{L}^2}\right) \leq c_8\|u\|_{\mathbb{H}^3},\label{ineq4}
\end{eqnarray}
where $c_1,\cdots,c_8>0$.

\section{Output Feedback Stabilization of the Linear KdV Equation}

To stabilizes \eqref{main}-\eqref{bcmain3}, we use the design presented in \cite{Cerpa1} for the linear system. Consider the following linear KdV equation
\begin{eqnarray}
u_t(x,t)+u_x(x,t)+u_{xxx}(x,t)&=&0\label{lkdv}
\end{eqnarray}
with the following boundary conditions
\begin{eqnarray}
u(0,t)&=&U(t),\;\;u(1,t)=0,\;\;u_x(1,t)=0\label{bclkdv3}
\end{eqnarray}
We assume only $u_{xx}(1,t)$ is measurable. If we select the control law $U(t)$ as
\begin{eqnarray}
U(t) &=& \int_0^1\!k(0,y)\hat{u}(y,t)\,\mathrm{d}y\label{claw}
\end{eqnarray}
where $\hat{u}$ is computed from
\begin{eqnarray}
\hat{u}_t+\hat{u}_x+\hat{u}_{xxx}+p_1(x)[y(t)-\hat{u}_{xx}(1,t)]&=&0\label{obslkdv}
\end{eqnarray}
and where $p_1(x)=p(x,1)$ with boundary conditions
\begin{eqnarray}
\hat{u}(0,t)&=&U(t),\;\;\hat{u}(1,t)=0,\;\;\hat{u}_x(1,t)=0\label{obsbclkdv3}
\end{eqnarray}
it can be shown that the origin of \eqref{lkdv}-\eqref{bclkdv3} is exponentially stable, where the kernels $k$ in \eqref{claw} is solution of the following kernel equation
\begin{eqnarray}
k_{xxx}+k_{yyy}+k_x+k_y &=& -\lambda k\label{ker1}
\end{eqnarray}
with boundary conditions
\begin{eqnarray}
k(x,1) &=& 0,\;\;k(x,x) = 0,\;\;k_x(x,x) = \frac{\lambda}{3}(1-x)\label{ker2}
\end{eqnarray}
Similarly, the kernel $p$ in \eqref{obslkdv} is solution of the following kernel equation
\begin{eqnarray}
p_{xxx}+p_{yyy}+p_x+p_y &=& \lambda p\label{ker3}
\end{eqnarray}
with boundary conditions
\begin{eqnarray}
p(x,x) &=& 0,\;\;p_x(x,x) = \frac{\lambda}{3}x,\;\;p(0,y) = 0\label{ker4}
\end{eqnarray}
Both kernel equations evolve in a triangular domain $\mathcal{T}=\left\{(x,y)|0\leq y\leq x\leq1\right\}$. The output-feedback boundary stabilization for the linear KdV equation is given as follow.
\begin{theorem}\citep{Cerpa1}
Consider systems \eqref{lkdv}-\eqref{bclkdv3} and \eqref{obslkdv}-\eqref{obsbclkdv3} with control law \eqref{claw} and initial conditions $u_0\in\mathbb{H}^3([0,1])$ and $\hat{u}_0\in\mathbb{L}^2([0,1])$. Then, there exists $\lambda>0$ and $c>0$ such that
\begin{eqnarray}
\|u(\cdot,t)\|_{\mathbb{H}^3}+\|\hat{u}(\cdot,t)\|_{\mathbb{L}^2} &\leq& ce^{-\lambda t}\left(\|u_0\|_{\mathbb{H}^3}+\|\hat{u}_0\|_{\mathbb{L}^2}\right)
\end{eqnarray}
\end{theorem}
\begin{remark}
We slightly modified the proof of this theorem so it can be used to analyze the stability of the KdV equation in section 4. The contribution (novelty) of this paper is in introduction of a strict Lyapunov functional for the output-feedback problem. To ease the reader, we present the proof of the above theorem as follow.
\end{remark}
\begin{proof}
Let the observer error be given by $\tilde{u}=u-\hat{u}$. We define new target variables $\hat{\omega}$ and $\tilde{\omega}$ using the following transformations
\begin{eqnarray}
\hat{\omega}(x,t) = \hat{u}(x,t) - \int_x^1\!k(x,y)\hat{u}(y,t)\,\mathrm{d}y\label{trans1}\\
\tilde{u}(x,t) = \tilde{\omega}(x,t) - \int_x^1\!p(x,y)\tilde{\omega}(y,t)\,\mathrm{d}y\label{trans2}
\end{eqnarray}
It can be shown that, if the kernels verify \eqref{ker1}-\eqref{ker4}, then $\hat{\omega}$ and $\tilde{\omega}$ satisfy the following equations
\begin{eqnarray}
\hat{\omega}_t+\hat{\omega}_x+\hat{\omega}_{xxx}+\lambda\hat{\omega} &=& -\bar{p}(x)\tilde{\omega}_{xx}(1,t)\label{miho1}\\
\hat{\omega}(0,t) = 0,\;\;\hat{\omega}(1,t) = 0,\;\;\hat{\omega}_x(1,t) &=& 0\\
\tilde{\omega}_t+\tilde{\omega}_x+\tilde{\omega}_{xxx}+\lambda\tilde{\omega} &=& 0\label{meh}\\
\tilde{\omega}(0,t) = 0,\;\;\tilde{\omega}(1,t) = 0,\;\;\tilde{\omega}_x(1,t) &=& 0\label{miho4}
\end{eqnarray}
where
\begin{eqnarray}
\bar{p}(x) = p_1(x)-\int_x^1\!k(x,y)p_1(y)\,\mathrm{d}y
\end{eqnarray}
Remark that, by bounding the norms in \eqref{meh}, we have
\begin{eqnarray}
\|\tilde{\omega}_t\|_{\mathbb{L}^2}&\leq&K_1\|\tilde{\omega}_{xxx}\|_{\mathbb{L}^2}\label{bing1}\\
\|\tilde{\omega}_{xxx}\|_{\mathbb{L}^2}&\leq&K_2\|\tilde{\omega}_t\|_{\mathbb{L}^2}\label{bing2}
\end{eqnarray}
where $K_1,K_2>0$. Thus, the $\|\tilde{\omega}_t\|_{\mathbb{L}^2}$ is equivalent to $\|\tilde{\omega}\|_{\mathbb{H}^3}$. Let us consider the following Lyapunov functional
\begin{eqnarray}
V(t) &=& \frac{A}{2}\int_0^1\!\hat{\omega}^2\,\mathrm{d}x+\frac{B}{2}\int_0^1\!\tilde{\omega}^2\,\mathrm{d}x+\frac{B}{2}\int_0^1\!\tilde{\omega}_t^2\,\mathrm{d}x
\end{eqnarray}
where $A$ and $B$ are positive constants. We can observe that the last term of the Lyapunov functional is equivalent to $\|\tilde{\omega}\|_{\mathbb{H}^3}$. Calculating the derivative of the Lyapunov functional along \eqref{miho1}-\eqref{miho4}, we have
\begin{eqnarray}
V(t) &\leq& A\left(-\lambda+\frac{D^2}{A}\right)\int_0^1\!\hat{\omega}^2\,\mathrm{d}x+A^2\tilde{\omega}_{xx}^2(1,t)\nonumber\\
&&-\lambda B\int_0^1\!\tilde{\omega}^2\,\mathrm{d}x-\lambda B\int_0^1\!\tilde{\omega}_t^2\,\mathrm{d}x
\end{eqnarray}
where $D=\max_{x\in[0,1]}\left\{p_1(x)-\int_x^1\!k(x,y)p_1(y)\,\mathrm{d}y\right\}$. From \eqref{bing2} and using integration by parts, we compute the bound of $\tilde{\omega}_{xx}^2(1,t)$ as follow
\begin{eqnarray}
|\tilde{\omega}_{xx}^2(1,t)| &\leq& a\|\tilde{\omega}\|^2_{\mathbb{L}^2}+b\|\tilde{\omega}_t\|^2_{\mathbb{L}^2}
\end{eqnarray}
Thus, we have
\begin{eqnarray}
\dot{V}(t) &\leq& A\left(-\lambda+\frac{D^2}{A}\right)\int_0^1\!\hat{\omega}^2\,\mathrm{d}x\nonumber\\
&&+B\left(-\lambda+\frac{aA^2}{B}\right)\int_0^1\!\tilde{\omega}^2\,\mathrm{d}x\nonumber\\
&&+B\left(-\lambda+\frac{bA^2}{B}\right)\int_0^1\!\tilde{\omega}_t^2\,\mathrm{d}x
\end{eqnarray}
Since we can choose arbitrary large $\lambda$, there exists $\epsilon>0$ such that
\begin{eqnarray}
\dot{V}(t) &\leq& -\epsilon V(t)
\end{eqnarray}
This completes the proof.
\end{proof}
The inverse of the transformations \eqref{trans1} and \eqref{trans2} are given by
\begin{eqnarray}
\hat{u}(x,t) &=& \hat{\omega}(x,t) + \int_x^1\!l(x,y)\hat{\omega}(y,t)\,\mathrm{d}y\label{invtrans1}\\
\tilde{\omega}(x,t) &=& \tilde{u}(x,t) + \int_x^1\!r(x,y)\tilde{u}(y,t)\,\mathrm{d}y\label{invtrans2}
\end{eqnarray}
The kernel $l(x,y)$ satisfy
\begin{eqnarray}
l_{xxx}+l_{yyy}+l_x+l_y &=& \lambda l\\
l(x,1) = 0,\;\;l(x,x) = 0,\;\;l_x(x,x) &=& \frac{\lambda}{3}(1-x)
\end{eqnarray}
while the kernel $r$ satisfy
\begin{eqnarray}
r_{xxx}+r_{yyy}+r_x+r_y &=& -\lambda r\label{kerr1}\\
r(x,x) = 0,\;\;r_x(x,x) = \frac{\lambda}{3}x,\;\;r(0,y) &=& 0\label{kerr2}
\end{eqnarray}
The kernel $k(x,y)$ and $l(x,y)$ are related by the formula
\begin{eqnarray}
l(x,y)-k(x,y) &=& \int_x^y\! k(x,\xi)l(\xi,y)\,\mathrm{d}\xi
\end{eqnarray}
A similar relation is also found in $p(x,y)$ and $r(x,y)$. The existence and uniqueness of the kernel solutions can be proved using the method of successive approximations.

\section{Output Feedback Stabilization of the KdV Equation}

We will show that the linear design \eqref{lkdv}-\eqref{bclkdv3} works for the nonlinear system \eqref{main}. Therefore, we design the following nonlinear observer as follow
\begin{eqnarray}
\hat{u}_t+\hat{u}_x+\hat{u}_{xxx}+\hat{u}\hat{u}_x+p_1(x)\tilde{u}_{xx}(1,t)&=&0\label{nonobs}
\end{eqnarray}
with boundary conditions
\begin{eqnarray}
\hat{u}(0,t)&=&U(t),\;\;\hat{u}(1,t)=0,\;\;\hat{u}_x(1,t)=0\label{bcnonobs3}
\end{eqnarray}
The observer error system is given by
\begin{eqnarray}
\tilde{u}_t+\tilde{u}_x+\tilde{u}_{xxx}+uu_x-\hat{u}\hat{u}_x-p_1(x)\tilde{u}_{xx}(1,t)&=&0\label{obssys}
\end{eqnarray}
with boundary conditions
\begin{eqnarray}
\tilde{u}(0,t)&=&0,\;\;\tilde{u}(1,t)=0,\;\;\tilde{u}_x(1,t)=0\label{obssysbc}
\end{eqnarray}
Let us define the following functionals
\begin{eqnarray}
\mathcal{K}[\omega] &=& \omega(x,t) - \int_x^1\!k(x,y)\omega(y,t)\,\mathrm{d}y\\
\mathcal{L}[\omega] &=& \omega(x,t) + \int_x^1\!l(x,y)\omega(y,t)\,\mathrm{d}y\\
\mathcal{P}[\omega] &=& \omega(x,t) - \int_x^1\!p(x,y)\omega(y,t)\,\mathrm{d}y\\
\mathcal{R}[\omega] &=& \omega(x,t) + \int_x^1\!r(x,y)\omega(x,t)\,\mathrm{d}y\\
\mathcal{L}_1[\omega] &=& l(x,x)\omega(x,t) + \int_x^1\!l_x(x,y)\omega(y,t)\,\mathrm{d}y\\
\mathcal{P}_1[\omega] &=& p(x,x)\omega(x,t)-\int_x^1\!p_x(x,y)\omega(y,t)\,\mathrm{d}y
\end{eqnarray}
By direct observation, these functionals satisfy
\begin{eqnarray}
|\mathcal{K}[\omega]| &\leq& C_1\left(|\omega|+\|\omega\|_{\mathbb{L}^1}\right)\\
|\mathcal{L}[\omega]| &\leq& C_2\left(|\omega|+\|\omega\|_{\mathbb{L}^1}\right)\\
|\mathcal{P}[\omega]| &\leq& C_3\left(|\omega|+\|\omega\|_{\mathbb{L}^1}\right)\\
|\mathcal{R}[\omega]| &\leq& C_4\left(|\omega|+\|\omega\|_{\mathbb{L}^1}\right)\\
|\mathcal{L}_1[\omega]| &\leq& C_5\left(|\omega|+\|\omega\|_{\mathbb{L}^1}\right)\\
|\mathcal{P}_1[\omega]| &\leq& C_6\left(|\omega|+\|\omega\|_{\mathbb{L}^1}\right)
\end{eqnarray}
for $C_1,\cdots,C_6>0$. Furthermore, we calculate the derivatives of \eqref{trans1} with respect to $x$
\begin{eqnarray}
\hat{\omega}_x(x,t) &=& \hat{u}_x(x,t) + k(x,x)\hat{u}(x,t)\nonumber\\
&&- \int_x^1\!k_x(x,y)\hat{u}(y,t)\,\mathrm{d}y\\
\hat{\omega}_{xx}(x,t) &=& \hat{u}_{xx}(x,t) + \frac{\mathrm{d}}{\mathrm{d}x}k(x,x)\hat{u}(x,t)+ k(x,x)\hat{u}_x(x,t) \nonumber\\
&&+k_x(x,x)\hat{u}(x,t)- \int_x^1\!k_{xx}(x,y)\hat{u}(y,t)\,\mathrm{d}y\\
\hat{\omega}_{xxx}(x,t) &=& \hat{u}_{xxx}(x,t)+ \frac{\mathrm{d}^2}{\mathrm{d}x^2}k(x,x)\hat{u}(x,t)\nonumber\\
&&+ 2\frac{\mathrm{d}}{\mathrm{d}x}k(x,x)\hat{u}_x(x,t)+ k(x,x)\hat{u}_{xx}(x,t) \nonumber\\
&&+\frac{\mathrm{d}}{\mathrm{d}x}k_x(x,x)\hat{u}(x,t)+k_x(x,x)\hat{u}_x(x,t)\nonumber\\
&&+k_{xx}(x,x)\hat{u}(x,t)- \int_x^1\!k_{xxx}(x,y)\hat{u}(y,t)\,\mathrm{d}y
\end{eqnarray}
and the derivative of \eqref{trans1} with respect to $t$ along \eqref{nonobs}
\begin{eqnarray}
\hat{\omega}_t(x,t) &=& -\hat{u}_x(x,t)-\hat{u}_{xxx}(x,t)+\int_x^1\!k(x,y)\hat{u}_y(y,t)\,\mathrm{d}y\nonumber\\
&&+k(x,1)\hat{u}_{xx}(1,t)-k(x,x)\hat{u}_{xx}(x,t)\nonumber\\
&&-k_y(x,1)\hat{u}_x(1,t)+k_y(x,x)\hat{u}_x(x,t)\nonumber\\
&&+k_{yy}(x,1)\hat{u}(1,t)-k_{yy}(x,x)\hat{u}(x,t)\nonumber\\
&&-\int_x^1\!k_{yyy}(x,y)\hat{u}(y,t)\,\mathrm{d}y\nonumber\\
&&-p_1(x)\tilde{u}_{xx}(1,t)+\int_x^1\!k(x,y)p_1(y)\,\mathrm{d}y\tilde{u}_{xx}(1,t)\nonumber\\
&&-\hat{u}\hat{u}_x(x,t)+\int_x^1\!k(x,y)\hat{u}(y,t)\hat{u}_y(y,t)\,\mathrm{d}y
\end{eqnarray}
Plugging the kernel \eqref{ker1}-\eqref{ker2} into the above equations, the observe $\hat{\omega}$ satisfy
\begin{eqnarray}
\hat{\omega}_t+\hat{\omega}_x+\hat{\omega}_{xxx}+\lambda\hat{\omega} &=& -\bar{p}(x)\tilde{\omega}_{xx}(1,t)-F[\hat{\omega},\hat{\omega}_x]\label{king1}
\end{eqnarray}
with boundary conditions
\begin{eqnarray}
\hat{\omega}(0,t) = 0,\;\;\hat{\omega}(1,t) &=& 0,\;\;\hat{\omega}_x(1,t) = 0\label{king2}
\end{eqnarray}
where
\begin{eqnarray}
F[\hat{\omega},\hat{\omega}_x] &=& \mathcal{K}[\mathcal{L}[\hat{\omega}]\left(\hat{\omega}_x+\mathcal{L}_1[\hat{\omega}]\right)]
\end{eqnarray}
This functional satisfy
\begin{eqnarray}
|F|&\leq& c_1\left(\|\hat{\omega}\|_{\mathbb{L}^2}+|\hat{\omega}|\right)\left(\|\hat{\omega}_x\|_{\mathbb{L}^2}+|\hat{\omega}_x|\right)\nonumber\\
&&+c_2\left(\|\hat{\omega}\|_{\mathbb{L}^2}^2+|\hat{\omega}|^2\right)
\end{eqnarray}
Next, computing the derivatives of \eqref{invtrans2} with respect to $x$, we have
\begin{eqnarray}
\tilde{\omega}_x(x,t) &=& \tilde{u}_x(x,t)-r(x,x)\tilde{u}(x,t) + \int_x^1\!r_x(x,y)\tilde{u}(y,t)\,\mathrm{d}y\nonumber\\
\\
\tilde{\omega}_{xx}(x,t) &=& \tilde{u}_{xx}(x,t)-\frac{\mathrm{d}}{\mathrm{d}x}r(x,x)\tilde{u}(x,t)-r(x,x)\tilde{u}_x(x,t) \nonumber\\
&&-r_x(x,x)\tilde{u}(x,t)+\int_x^1\!r_{xx}(x,y)\tilde{u}(y,t)\,\mathrm{d}y\\
\tilde{\omega}_{xxx}(x,t) &=& \tilde{u}_{xxx}(x,t)-\frac{\mathrm{d}^2}{\mathrm{d}x^2}r(x,x)\tilde{u}(x,t)\nonumber\\
&&-2\frac{\mathrm{d}}{\mathrm{d}x}r(x,x)\tilde{u}_x(x,t)-r(x,x)\tilde{u}_{xx}(x,t) \nonumber\\
&&-\frac{\mathrm{d}}{\mathrm{d}x}r_x(x,x)\tilde{u}(x,t)-r_x(x,x)\tilde{u}_x(x,t)\nonumber\\
&&-r_{xx}(x,y)\tilde{u}(y,t)+\int_x^1\!r_{xxx}(x,y)\tilde{u}(y,t)\,\mathrm{d}y
\end{eqnarray}
Furthermore, we calculate
\begin{eqnarray}
\tilde{\omega}_t(x,t) &=&-\tilde{u}_x(x,t)-\tilde{u}_{xxx}(x,t)+p_1(x)\tilde{u}_{xx}(1,t) \nonumber\\
&&-r(x,1)\tilde{u}(1,t)+r(x,x)\tilde{u}(x,t)\nonumber\\
&&-r(x,1)\tilde{u}_{yy}(1,t)+r(x,x)\tilde{u}_{yy}(x,t)\nonumber\\
&&+r_y(x,1)\tilde{u}_{y}(1,t)-r_y(x,x)\tilde{u}_{y}(x,t)\nonumber\\
&&-r_{yy}(x,1)\tilde{u}(1,t)+r_{yy}(x,x)\tilde{u}(x,t)\nonumber\\
&&+\int_x^1\!\left(r_y(x,y)+r_{yyy}(x,y)\right)\tilde{u}(y,t)\,\mathrm{d}y\nonumber\\
&&+\int_x^1\!r(x,y)p_1(y)\,\mathrm{d}y\tilde{u}_{yy}(1,t)\nonumber\\
&&-uu_x+\hat{u}\hat{u}_x+\int_x^1\!r(x,y)\left(-uu_y+\hat{u}\hat{u}_y\right)\,\mathrm{d}y
\end{eqnarray}
Plugging the kernel equation \eqref{kerr1}-\eqref{kerr2}, we have
\begin{eqnarray}
\tilde{\omega}_t+\tilde{\omega}_x+\tilde{\omega}_{xxx}+\lambda\tilde{\omega} &=& -G[\hat{\omega},\hat{\omega}_x,\tilde{\omega},\tilde{\omega}_x]\label{king3}
\end{eqnarray}
with boundary conditions
\begin{eqnarray}
\tilde{\omega}(0,t) = 0,\;\;\tilde{\omega}(1,t) &=& 0,\;\;\tilde{\omega}_x(1,t) = 0\label{king4}
\end{eqnarray}
where
\begin{eqnarray}
G[\hat{\omega},\hat{\omega}_x,\tilde{\omega},\tilde{\omega}_x] &=& \mathcal{R}[\mathcal{P}[\tilde{\omega}]\left(\tilde{\omega}_x+\mathcal{P}_1[\tilde{\omega}]\right)]\nonumber\\
&&+\mathcal{R}[\mathcal{P}[\tilde{\omega}]\left(\hat{\omega}_x+\mathcal{L}_1[\hat{\omega}]\right)]\nonumber\\
&&+\mathcal{R}[\mathcal{L}[\hat{\omega}]\left(\tilde{\omega}_x+\mathcal{P}_1[\tilde{\omega}]\right)]\nonumber\\
&&-\mathcal{R}[\mathcal{L}[\hat{\omega}]\left(\hat{\omega}_x+\mathcal{L}_1[\hat{\omega}]\right)]
\end{eqnarray}
This functional satisfy
\begin{eqnarray}
|G| &\leq& c_1\left(\|\tilde{\omega}\|_{\mathbb{L}^2}+|\tilde{\omega}|\right)\left(\|\tilde{\omega}_x\|_{\mathbb{L}^2}+|\tilde{\omega}_x|\right)\nonumber\\
&&+c_2\left(\|\tilde{\omega}\|_{\mathbb{L}^2}^2+|\tilde{\omega}|^2\right)+c_3\left(\|\tilde{\omega}_x\|_{\mathbb{L}^2}^2+|\tilde{\omega}_x|^2\right)\nonumber\\
&&+c_4\left(\|\hat{\omega}\|_{\mathbb{L}^2}+|\hat{\omega}|\right)\left(\|\hat{\omega}_x\|_{\mathbb{L}^2}+|\hat{\omega}_x|\right)\nonumber\\
&&+c_5\left(\|\hat{\omega}\|_{\mathbb{L}^2}^2+|\hat{\omega}|^2\right)+c_6\left(\|\hat{\omega}_x\|_{\mathbb{L}^2}^2+|\hat{\omega}_x|^2\right)
\end{eqnarray}
We can observe that \eqref{king1}-\eqref{king4} are PDE-PDE cascade systems. To study its stability, first we denote $\eta=\omega_t$. Taking a derivative of \eqref{king1}-\eqref{king4} with respect to $t$, we have
\begin{eqnarray}
\hat{\eta}_t+\hat{\eta}_x+\hat{\eta}_{xxx}+\lambda\hat{\eta} &=& -\bar{p}(x)\tilde{\eta}_{xx}(1,t)-\mathcal{L}[\hat{\eta}]\hat{\omega}_x\nonumber\\
&&-F_1[\hat{\omega},\hat{\eta},\hat{\eta}_x]\label{kong1}\\
\hat{\eta}(0,t) = 0,\;\;\hat{\eta}(1,t) &=& 0,\;\;\hat{\eta}_x(1,t) = 0\\
\tilde{\eta}_t+\tilde{\eta}_x+\tilde{\eta}_{xxx}+\lambda\tilde{\eta} &=& -\mathcal{P}[\tilde{\eta}]\tilde{\omega}_x-\mathcal{P}[\tilde{\eta}]\hat{\omega}_x-\mathcal{L}[\hat{\eta}]\tilde{\omega}_x\\
&&+\mathcal{L}[\hat{\eta}]\hat{\omega}_x-G_1[\hat{\omega}, \hat{\eta},\hat{\eta}_x,\tilde{\omega},\tilde{\eta},\tilde{\eta}_x]\nonumber\label{kong3}\\
\tilde{\eta}(0,t) = 0,\;\;\tilde{\eta}(1,t) &=& 0,\;\;\tilde{\eta}_x(1,t) = 0\label{kong4}
\end{eqnarray}
where
\begin{eqnarray}
F_1[\hat{\omega},\hat{\eta},\hat{\eta}_x] &=& -\mathcal{L}[\hat{\eta}(1,t)]\hat{\omega}(1,t)+\mathcal{L}[\hat{\eta}(x,t)]\hat{\omega}(x,t)\nonumber\\
&&+\int_x^1\!\left(\omega_x+\mathcal{L}_1[\hat{\eta}]\right)\hat{\omega}\,\mathrm{d}y+\mathcal{K}[\mathcal{L}[\hat{\eta}]\mathcal{L}_1[\hat{\omega}]]\nonumber\\
&&+\mathcal{K}[\mathcal{L}[\hat{\omega}]\left(\hat{\eta}_x+\mathcal{L}_1[\hat{\eta}]\right)]
\end{eqnarray}
and
\begin{eqnarray}
&&G_1[\hat{\omega}, \hat{\eta},\hat{\eta}_x,\tilde{\omega},\tilde{\eta},\tilde{\eta}_x] \nonumber\\
&=&-\mathcal{P}[\tilde{\eta}(x,t)]\tilde{\omega}(x,t)-\int_x^1\!\mathcal{P}_x[\tilde{\eta}]\tilde{\omega}\,\mathrm{d}y\nonumber\\
&&+\mathcal{R}[\mathcal{P}[\tilde{\eta}]\mathcal{P}_1[\tilde{\omega}]]+\mathcal{R}[\mathcal{P}[\tilde{\omega}]\left(\tilde{\eta}_x+\mathcal{P}_1[\tilde{\eta}]\right)]\nonumber\\
&&-\mathcal{P}[\tilde{\eta}(x,t)]\hat{\omega}(x,t)-\int_x^1\!\mathcal{P}_x[\tilde{\eta}]\hat{\omega}\,\mathrm{d}y\nonumber\\
&&+\mathcal{R}[\mathcal{P}[\tilde{\eta}]\mathcal{L}_1[\hat{\omega}]]+\mathcal{R}[\mathcal{P}[\tilde{\omega}]\left(\hat{\eta}_x+\mathcal{L}_1[\hat{\eta}]\right)]\nonumber\\
&&-\mathcal{L}[\hat{\eta}(x,t)]\tilde{\omega}(x,t)-\int_x^1\!\mathcal{L}_x[\hat{\eta}]\tilde{\omega}\,\mathrm{d}y\nonumber\\
&&+\mathcal{R}[\mathcal{L}[\hat{\eta}]\mathcal{P}_1[\tilde{\omega}]]+\mathcal{R}[\mathcal{L}[\hat{\omega}]\left(\tilde{\eta}_x+\mathcal{P}_1[\tilde{\eta}]\right)]\nonumber\\
&&+\mathcal{L}[\hat{\eta}(x,t)]\hat{\omega}(x,t)+\int_x^1\!\mathcal{L}_x[\hat{\eta}]\hat{\omega}\,\mathrm{d}y\nonumber\\
&&-\mathcal{R}[\mathcal{L}[\hat{\eta}]\mathcal{L}_1[\hat{\omega}]]-\mathcal{R}[\mathcal{L}[\hat{\omega}]\left(\hat{\eta}_x+\mathcal{L}_1[\hat{\eta}]\right)]
\end{eqnarray}
Bounding the norms for small $\|\hat{\omega}\|_{\infty}+\|\tilde{\omega}\|_{\infty}$ in \eqref{kong1} and \eqref{kong3}, we can prove that the norm $\|\hat{\omega}\|_{\mathbb{H}^3}+\|\tilde{\omega}\|_{\mathbb{H}^3}$ is equivalent to $\|\hat{\eta}\|_{\mathbb{L}^2}+\|\tilde{\eta}\|_{\mathbb{L}^2}=\|\hat{\omega}_t\|_{\mathbb{L}^2}+\|\tilde{\omega}_t\|_{\mathbb{L}^2}$. The main result of this paper is stated as follow.
\begin{theorem}
Consider systems \eqref{main}-\eqref{bcmain3} and \eqref{nonobs}-\eqref{bcnonobs3} with control law \eqref{claw} and initial conditions $u_0,\hat{u}_0\in\mathbb{H}^3([0,1])$. Then, there exists $\delta$, $\lambda>0$, and $c>0$ such that if $\|u_0\|_{\mathbb{H}^3}+\|\hat{u}_0\|_{\mathbb{H}^3}\leq\delta$, then
\begin{eqnarray}
\|u(\cdot,t)\|_{\mathbb{H}^3}+\|\hat{u}(\cdot,t)\|_{\mathbb{H}^3} &\leq& ce^{-\lambda t}\left(\|u_0\|_{\mathbb{H}^3}+\|\hat{u}_0\|_{\mathbb{H}^3}\right)
\end{eqnarray}
\end{theorem}

\begin{remark}
The difference between the linear system and the nonlinear system results lie in the smallness of the initial condition. Thus, in the nonlinear system, we only achieved local exponential stability.
\end{remark}

\begin{proof}
We introduce the following Lyapunov functional
\begin{eqnarray}
W(t) &=& \frac{A}{2}\int_0^1\!\hat{\omega}^2\,\mathrm{d}x+\frac{A}{2}\int_0^1\!\hat{\eta}^2\,\mathrm{d}x\nonumber\\
&&+\frac{B}{2}\int_0^1\!\tilde{\omega}^2\,\mathrm{d}x+\frac{B}{2}\int_0^1\!\tilde{\eta}^2\,\mathrm{d}x
\end{eqnarray}
Computing its derivative along \eqref{king1}-\eqref{king4} with respect to $t$, yields
\begin{eqnarray}
\dot{W}(t) &\leq& -\mu W(t)-A\int_0^1\!\hat{\omega}F\,\mathrm{d}x-A\int_0^1\!\hat{\eta}F_1\,\mathrm{d}x\\
&&-B\int_0^1\!\tilde{\omega}G\,\mathrm{d}x-B\int_0^1\!\tilde{\eta}G_1\,\mathrm{d}x-\int_0^1\!\hat{\eta}\mathcal{L}[\hat{\eta}]\hat{\omega}_x\,\mathrm{d}x\nonumber\\
&&-\int_0^1\!\tilde{\eta}\left(\mathcal{P}[\tilde{\eta}]\tilde{\omega}_x+\mathcal{P}[\tilde{\eta}]\hat{\omega}_x+\mathcal{L}[\hat{\eta}]\tilde{\omega}_x-\mathcal{L}[\hat{\eta}]\hat{\omega}_x\right)\,\mathrm{d}x\nonumber
\end{eqnarray}
We estimate
\begin{eqnarray}
&&\left|-\int_0^1\!\hat{\eta}\mathcal{L}[\hat{\eta}]\hat{\omega}_x\,\mathrm{d}x\right|\leq K_1\|\hat{\eta}\|_{\infty}W(t)\\
&&\left|-\int_0^1\!\tilde{\eta}\left(\mathcal{P}[\tilde{\eta}]\tilde{\omega}_x+\mathcal{P}[\tilde{\eta}]\hat{\omega}_x+\mathcal{L}[\hat{\eta}]\tilde{\omega}_x-\mathcal{L}[\hat{\eta}]\hat{\omega}_x\right)\,\mathrm{d}x\right| \leq\nonumber\\
&&K_2\|\tilde{\eta}\|_{\infty}W(t)
\end{eqnarray}
Furthermore, we estimate
\begin{eqnarray}
\left|-A\int_0^1\!\hat{\omega}F\,\mathrm{d}x\right|&\leq&K_3\|\hat{\omega}_x\|_{\mathbb{\infty}}W(t)\\
\left|-A\int_0^1\!\hat{\eta}F_1\,\mathrm{d}x\right|&\leq&K_4\left(\|\hat{\eta}\|W(t)+W(t)^{3/2}\right)\\
\left|-B\int_0^1\!\tilde{\omega}G\,\mathrm{d}x\right|&\leq&K_5\|\tilde{\omega}_x\|_{\mathbb{\infty}}W(t)\nonumber\\
\left|-B\int_0^1\!\tilde{\eta}G_1\,\mathrm{d}x\right|&\leq&K_6\left(\|\tilde{\eta}\|W(t)+W(t)^{3/2}\right)
\end{eqnarray}
where $K_1,\cdots,K_6>0$. Thus, we have
\begin{eqnarray}
\dot{W}(t) &\leq& -\mu W(t)+CW(t)^{\frac{3}{2}}
\end{eqnarray}
for some positive $\mu$ and $C$. Then, for any $\mu_0$ such that $0<\mu_0<\mu$, there exists $\delta_0$ such that
\begin{eqnarray}
C\left|W^{3/2}\right| < \left(\mu-\mu_0\right)W,\;\forall W<\delta_0,
\end{eqnarray}
which implies that
\begin{eqnarray}
\dot{W}<-\mu_0W,\;\forall W<\delta_0.
\end{eqnarray}
Since $W$ is equivalent to $\|\omega\|_{\mathbb{H}^3}+\|\hat{\omega}\|_{\mathbb{H}^3}$ when $\|\omega\|_{\mathbb{\infty}}+\|\hat{\omega}\|_{\mathbb{\infty}}$ is sufficiently small, this concludes the proof.
\end{proof}

\section{Numerical Example}

To show our linear control law works for the nonlinear system, we simulate the closed-loop system \eqref{main}-\eqref{bcmain3} and \eqref{nonobs}-\eqref{bcnonobs3} with control law \eqref{claw}. The initial conditions are chosen such that they are compatible with the boundary conditions. The result is presented in figure 1. We can observe in the controlled case, the controller drives the closed-loop system into its equilibrium.
\begin{figure}[h!]
  \centering
      \includegraphics[width=0.5\textwidth]{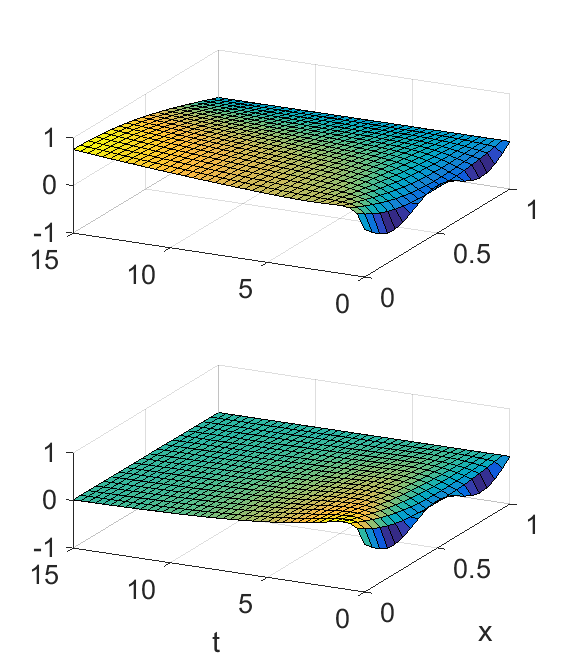}
  \caption{Uncontrolled case (top) and controlled case (bottom).}
\end{figure}

\section{Conclusion}

In this paper, we have presented output-feedback boundary stabilization of the KdV equation with actuation and measurement on only one boundary. The control law was obtained using backstepping method for the linear system. Using a strict Lyapunov functional, we have shown local exponential $\mathbb{H}^3$ stability of the state and of the observer error.

\balance
\bibliography{ifacconf}

\end{document}